\documentclass{amsart}%
\usepackage{amsmath}
\usepackage{amsfonts}
\usepackage{amssymb}
\usepackage{amscd}
\usepackage[dvipsnames]{xcolor}
\usepackage{graphicx, enumerate}%
\usepackage[justification=centering]{caption}
\usepackage[nameinlink,noabbrev]{cleveref}
\usepackage{caption}
\usepackage{subcaption}
\usepackage{multicol}
\usepackage{float}
\usepackage{stackengine}

\providecommand{\U}[1]{\protect \rule{.1in}{.1in}}
\newtheorem{thm}{Theorem}[section]
\newtheorem{lem}[thm]{Lemma}
\newtheorem{prop}[thm]{Proposition}

\numberwithin{equation}{section}

\theoremstyle{definition}
\newtheorem{question}{Question}

\let\max\relax \DeclareMathOperator*\max{\vphantom{p}max}
\newcommand{\tr}{\operatorname{tr}}

\newcommand{\Int}{\operatorname{Int}}

\title{On a question of Erd\H{o}s on doubly stochastic matrices}

\author{Ludovick Bouthat}
 \address{D\'epartement de math\'ematiques et de statistique, Universit\'e Laval, 1045, avenue de la M\'edecine, Qu\'ebec, QC, G1V\;0A6, Canada}
 \email{Ludovick.Bouthat.1@ulaval.ca}

 \author{Javad Mashreghi}
 \address{D\'epartement de math\'ematiques et de statistique, Universit\'e Laval, 1045, avenue de la M\'edecine, Qu\'ebec, QC, G1V\;0A6, Canada}
 \email{Javad.Mashreghi@mat.ulaval.ca}

 \author{Fr\'ed\'eric Morneau-Gu\'erin}
 \address{D\'epartement \'Education, Universit\'e T\'ELUQ, Qu\'ebec, QC, G1K 9H6, Canada.}
 \email{Frederic.Morneau-Guerin@teluq.ca}



\begin{document}

\begin{abstract}
In a celebrated paper of Marcus and Ree (1959), it was shown that if  $A=[a_{ij}]$ is an $n \times n$ doubly stochastic matrix, then there is a permutation $\sigma \in S_n$ such that $\sum_{i,j=1}^{n} a_{i,j}^{2} \leq \sum_{i=1}^{n} a_{i,\sigma(i)}$.  Erd\H{o}s asked for which doubly stochastic matrices the inequality is saturated. Although Marcus and Ree provided some insight for the set of solutions, the question appears to have fallen into oblivion. Our goal is to provide a complete answer in the particular, yet non-trivial, case when $n=3$.
\end{abstract}

\maketitle

%%%%%%%%%%%%%%%%%%%%%%%%%%%%%%%%%%%%%%%%%%%%%%%%%%%%%%%%%%%%%%%
%%%%%%%%%%%%%%%%%%%%%%%%%%%%%%%%%%%%%%%%%%%%%%%%%%%%%%%%%%%%%%%
%%%%%%%%%%%%%%%%%%%%%%%%%%%%%%%%%%%%%%%%%%%%%%%%%%%%%%%%%%%%%%%
%%%%%%%%%%%%%%%%%%%%%%%%%%%%%%%%%%%%%%%%%%%%%%%%%%%%%%%%%%%%%%%

\section{Introduction}

A square matrix $A= [a_{i,j}]$ is called \textit{doubly stochastic} if every entry $a_{i,j}$ is non-negative and all row sums and all column sums of $A$ are equal to 1, i.e.,
$$
a_{i,j} \geq 0, \quad\quad \sum\limits_{i=1}^n a_{i,j}= \sum\limits_{j=1}^n a_{i,j}=1,
$$
for all $i, j = 1, \dots, n$. Trivial examples of doubly stochastic matrices are the permutation matrices, in particular the identity matrix, and the matrix $J_n$ for which all entries are $\frac{1}{n}$. Let $\Omega_n$ be the set of all doubly stochastic matrices. The classic Birkhoff--von Neumann theorem says that $\Omega_n$ is the convex hull of the set of all $n \times n$ permutation matrices, and furthermore that the vertices of $\Omega_n$ are precisely the permutation matrices \cite{Birkhoff1946}.

Let us now make a brief digression. The \textit{permanent} of a square matrix is a polynomial function in the entries of the matrix similar to the determinant, but without the minus signs, i.e.,
$$
\text{perm}(A) \,:=~ \sum\limits_{\sigma \in S_n}\prod\limits_{i=1}^n a_{i,\sigma(i)}.
$$
It is straightforward to see there are $n!$ such permutations. Although the definition of the permanent of a given matrix differs from that of its determinant only in the signature of permutations, calculating permanents turns out to be computationally much more difficult than calculating determinants.

In 1926, van der Waerden \cite{van1926aufgabe}, via considerations outside the field of matrix analysis, proposed to determine the minimal permanent among all $n \times n$ doubly stochastic matrices. As $J_n$ is the unique trivial guess for the minima and that
\[
\text{perm}(J_n) \,=\, n! \, n^{-n},
\]
the assertion
$$
\big(A \in \Omega_n \text{ and } A \neq J_n\big) \,\Longrightarrow\, \big(\text{perm}(A) > n! \, n^{-n} \big)
$$
became known as the \textit{van der Waerden conjecture} (as it turned out, unbeknownst to van der Waerden himself for several decades \cite{MR672919}).

The theory of permanents was a very active field of research from 1950s to 1980s, with the van der Waerden conjecture  responsible for the bulk of the research papers. After mystifying mathematicians for over half a century, the van der Waerden's conjecture ultimately gave in -- as so often happens in mathematics -- from multiple directions at the same time. The first proof of the conjecture to appear in print (in Russian at the end of 1980 \cite{MR602332} and in english the following year \cite{MR642395}) is the work of Egorychev.  Then, a few months later, to almost everyone's surprise, a paper submitted on May 14, 1979, by Falikman, with a completely different proof of the conjecture \cite{MR625097}. In 1958, Marcus and Ree  penned a seminal paper \cite{MR117243} that set off a large wave of articles on the van der Waerden conjecture \cite{MR104679}. In that paper, the two men studied properties of the diagonals of doubly stochastic matrices.

The line of research between Marcus and Ree was motivated by a fruitful discussion that they had with Erd\H{o}s about the van der Waerden  conjecture. It appears that, following the investigation carried out by Marcus and Ree, one of the questions asked by Erd\H{o}s faded into oblivion. The purpose of this article is to raise awareness on Erd\H{o}s' question, which is described in Section \ref{S:erdos-question} and its historical development explained in Section \ref{S:discussion-erdos-question}. Moreover, we completely answer this difficult question for the $3 \times 3$ doubly stochastic matrices. Our main result is stated in Section \ref{S:Main-Result}. Then, in Section \ref{Section3}, we work out a complete characterization of those $3 \times 3$ doubly stochastic matrices for which the Frobenius norm squared is equal to one of the diagonal sums (a weak form of Erd\H{o}s' question). Finally, in Section \ref{S:Prrof-Main-Result}, we identify the subset of this family of $3 \times 3$ doubly stochastic matrices for which this diagonal sum is precisely the maximal trace (Erd\H{o}s' question).

%%%%%%%%%%%%%%%%%%%%%%%%%%%%%%%%%%%%%%%%%%%%%%%%%%%%%%%%%%%%%%%
%%%%%%%%%%%%%%%%%%%%%%%%%%%%%%%%%%%%%%%%%%%%%%%%%%%%%%%%%%%%%%%
%%%%%%%%%%%%%%%%%%%%%%%%%%%%%%%%%%%%%%%%%%%%%%%%%%%%%%%%%%%%%%%
%%%%%%%%%%%%%%%%%%%%%%%%%%%%%%%%%%%%%%%%%%%%%%%%%%%%%%%%%%%%%%%

\section{The Erd\H{o}s question} \label{S:erdos-question}
At the time of collaboration between Marvin Marcus and Rimhak Ree on the van der Waerden conjecture, Erd\H{o}s pointed out that, given an $n \times n$ doubly stochastic matrix $A$, if the conjecture holds true (which was not yet an established fact back at the time), then the following statements hold.
\begin{enumerate}[(i)]
\item There is a permutation $\sigma$ such that
\begin{equation}\label{MarcusReeA}
\prod\limits_{i=1}^n a_{i,\sigma(i)} \geq n^{-n}.
\end{equation}
\item There is a permutation $\sigma$ such that
\begin{equation}\label{MarcusReeB}
\sum\limits_{i=1}^n a_{i, \sigma(i)} \geq 1 \quad \mbox{and} \quad a_{i,\sigma(i)} > 0 \quad \mbox{for} \quad i = 1, \dots, n.
\end{equation}
\item Any doubly stochastic matrix has a diagonal sum greater than or equal to 1.
\end{enumerate}

The last part needs clarification. Let $\sigma$ be any of the $n!$ permutations of the integers $1, 2, \dots, n$. We shall say that the sequence $a_{1, \sigma(1)},a_{2, \sigma(2)}, \dots, a_{n, \sigma(n)}$, consisting of exactly one element from each row and each column of $A$, is a \textit{diagonal} of $A$ corresponding to $\sigma$. Therefore, the third statement is a slightly weaker form of \eqref{MarcusReeB}. Marcus and Ree were unable to prove statement \eqref{MarcusReeA} (that is to say without invoking the van der Waerden's conjecture). However, they succeeded to prove the following considerably strengthened version of the third statement.

\begin{thm}[Marcus--Ree \cite{MR104679}]\label{ErdosMarcusRee}
If $A$ is an $n \times n$ doubly stochastic matrix, then there exists a permutation $\sigma$ such that
\[
\sum_{i=1}^n\sum_{j=1}^n a_{i,j}^2 \,\leq\, \sum_{i=1}^n a_{i,\sigma(i)}.
\]
\end{thm}
\noindent Recalling the Frobenius norm
\[
\|A\|_{\text{F}} \,=\, \left( \sum_{i=1}^n\sum_{j=1}^n a_{i,j}^2 \right)^{1/2},
\]
Marcus and Ree result says that any doubly stochastic matrix has a diagonal sum greater than or equal to its Frobenius norm squared. To easily see that this is indeed a stronger result than the statement \eqref{MarcusReeB}, one just need to observe that the Frobenius norm of any matrix is greater than or equal to its spectral radius which, in turn, is greater than or equal to 1 for all doubly stochastic matrices. In fact, it is straightforward to check that for any $A \in \Omega_n$, 1 is an eigenvalue corresponding to the eigenvector $x=(1,1,\dots,1)$.

For the purpose of our discussion, we proceed to encapsulate the essence of  \Cref{ErdosMarcusRee} in a definition. For an $n \times n$ matrix $A=[a_{i,j}]$, we define
\begin{equation}\label{E:max-tace}
\max\limits_{\tr}(A) \,=\, \max\limits_{\sigma \in S_n}\sum_{i=1}^n a_{i,\sigma(i)}.
\end{equation}
This quantity has drawn attention of leading scientists of the numerical and applied mathematics community over the years due to its relevance to the assignment problem in linear programming \cite{MR540958}. Note that Wang \cite{Wang1974} called \eqref{E:max-tace} the \textit{maximum diagonal sum} of $A$, but we preferred to call it the \emph{maximal trace} of $A$. Hence, Theorem \ref{ErdosMarcusRee} can be written as
\begin{equation}\label{E:EMR}
\|A\|_{\text{F}}^2 \,\leq\, \max\limits_{\tr}(A),
\end{equation}
where $A$ is a doubly stochastic matrix. Then, upon hearing about \Cref{ErdosMarcusRee},  Erd\H{o}s raised the following question:
\begin{question} (Erd\H{os})
When is $\|A\|_{\text{F}}^2 = \max\limits_{\tr}(A)$?
\end{question}
The question appears in the paper of Marcus and Ree \cite{MR104679}, along with some discussions, which we briefly present in Section \ref{S:discussion-erdos-question}. In the following, we say that $A$ \textit{saturates} the inequality \eqref{E:EMR} if $\|A\|_{\text{F}}^2  = \max\limits_{\tr}(A)$.

In our approach, in order to answer Erd\H{o}s' question, we need to treat primarily a weaker question in which we ask to characterize the doubly stochastic matrices $A$ for which there is a permutation $P$ such that
\[
\|A\|_{\text{F}}^2 \,=\, \tr(AP).
\]
However, note that such a matrix does not necessarily saturate the inequality \eqref{E:EMR}.

%%%%%%%%%%%%%%%%%%%%%%%%%%%%%%%%%%%%%%%%%%%%%%%%%%%%%%%%%%%%%%%
%%%%%%%%%%%%%%%%%%%%%%%%%%%%%%%%%%%%%%%%%%%%%%%%%%%%%%%%%%%%%%%
%%%%%%%%%%%%%%%%%%%%%%%%%%%%%%%%%%%%%%%%%%%%%%%%%%%%%%%%%%%%%%%
%%%%%%%%%%%%%%%%%%%%%%%%%%%%%%%%%%%%%%%%%%%%%%%%%%%%%%%%%%%%%%%

\section{Some preliminary discussions on Erd\H{o}s' question} \label{S:discussion-erdos-question}

One can easily answer Erd\H{o}s' question when $n = 2$. A $2 \times 2$ doubly stochastic matrix is of the following form
\[
A \,=\,
\begin{bmatrix}
t & 1-t \\
1-t & t
\end{bmatrix},
\]
where $0 \leq t \leq 1$.

\noindent Case I: It is straighforward to directly verify that the inequality \eqref{E:EMR} is saturated if $t \in \{0, 1/2, 1 \}$. These cases correspond to the two $2 \times 2$ permutation matrices
\[
\begin{bmatrix} 1 & 0 \\ 0 & 1 \end{bmatrix}
\quad \mbox{and} \quad
\begin{bmatrix} 0 & 1 \\ 1 & 0 \end{bmatrix},
\]
and to
\[
J_2 \,=\,
\begin{bmatrix}
\frac{1}{2} & \frac{1}{2} \\[3pt]
\frac{1}{2} & \frac{1}{2}
\end{bmatrix}.
\]

\noindent Case II: Assume that $0<t<1/2$. Then
\[
\|A\|_{\text{F}}^2 \,=\, 2t^2+2(1-t)^2,
\]
and
\[
a_{1,1}+a_{2,2} \,=\, 2t
\quad \mbox{and} \quad
a_{1,2}+a_{2,1} \,=\, 2(1-t),
\]
and thus
\[
\max\limits_{\tr}(A) \,=\, \max\{ 2t,\, 2(1-t) \} \,=\, 2(1-t).
\]
Hence,
\[
\max\limits_{\tr}(A)- \|A\|_{\text{F}}^2 \,=\, 2t(1-2t) \,>\, 0.
\]
\noindent Case III: Assume that $1/2<t<1$. Similar calculation as in the Case II.

For $n > 2$, characterizing the $n \times n$ doubly stochastic matrices that saturate the inequality \eqref{E:EMR} is no longer so simple. It is easy to see that for any permutation matrix $P$,
\[
\|P\|_{\text{F}}^2 \,=\, \max\limits_{\tr}(P) \,=\, n,
\]
and, as well, for the matrix
\[
J_n \,=\,
\begin{bmatrix}
\frac{1}{n} & \frac{1}{n} & \cdots & \frac{1}{n} \\[5pt]
\frac{1}{n} & \frac{1}{n} & \cdots & \frac{1}{n} \\
\vdots & \vdots & \ddots & \vdots \\
\frac{1}{n} & \frac{1}{n} & \cdots & \frac{1}{n} \\
\end{bmatrix},
\]
we have
\[
\|J_n\|_{\text{F}}^2 \,=\, \max\limits_{\tr}(J_n) \,=\, 1.
\]
In technical language, they saturate the inequality \eqref{E:EMR}. However, they are no longer the only $n \times n$ doubly stochastic matrices in doing so. Although Marcus and Ree did not provide a complete and definitive answer to Erd\H{o}s' question, they were nevertheless presented some interesting partial results. We mention two of them below. The first is a characterization of a family of doubly stochastic matrices that do not saturate the inequality \eqref{E:EMR}.

\begin{prop}[Marcus--Ree \cite{MR117243}, Corollary 2]\label{ErdosMarcusReePositive}
If $A$ is an $n \times n$ doubly stochastic matrix all of whose entries are strictly positive, and if $A \neq J_n$, then $A$ does not saturate the inequality \eqref{E:EMR}.
\end{prop}

Furthermore, Marcus and Ree also identified a family of $n \times n$ doubly stochastic matrices that saturate the inequality \eqref{E:EMR}. In the following, we denote by $\oplus$ the direct sum of two matrices, i.e., $A\oplus B = \left[ \begin{smallmatrix} A&0\\0&B \end{smallmatrix}\right]$.

\begin{prop}[Marcus--Ree \cite{MR117243}, Corollary 4]\label{ErdosMarcusReeThm}
If $A=P(J_{n_1} \oplus \cdots \oplus J_{n_r})Q$, where $P$ and $Q$ are permutation matrices and $n_1+\cdots+n_r=n$, then $A$ saturates the inequality \eqref{E:EMR}.
\end{prop}

In the light of Proposition \ref{ErdosMarcusReeThm}, it would be natural to expect that the doubly stochastic matrices saturating the inequality \eqref{E:EMR} are precisely those of the form $P(J_{n_1} \oplus \cdots \oplus J_{n_r})Q$, for some permutation matrices $P$ and $Q$ and some positive integers $n_1, \dots, n_r$ verifying $n_1+\cdots+n_r=n$. However, Marcus and Ree provided a counter-example to this hasty conjecture. The matrix
\begin{align}\label{eg - 1}
S \,=\,
\begin{bmatrix}
0 & \frac{1}{2} & \frac{1}{2} \\[3pt]
\frac{1}{2} & \frac{1}{4} & \frac{1}{4} \\[3pt]
\frac{1}{2} & \frac{1}{4} & \frac{1}{4}
\end{bmatrix}
\end{align}
saturates the inequality \eqref{E:EMR}, but it is  not of the form $P(J_{n_1} \oplus \cdots \oplus J_{n_r})Q$. It is interesting to note, however,  that it can be factored as the product of two matrices of the aforementioned form, i.e.,
\begin{equation} \label{E:matrix-prod}
\begin{bmatrix}
0 & \frac{1}{2} & \frac{1}{2} \\[3pt]
\frac{1}{2} & \frac{1}{4} & \frac{1}{4} \\[3pt]
\frac{1}{2} & \frac{1}{4} & \frac{1}{4}
\end{bmatrix}
\,=\,
\begin{bmatrix}
1 & 0 & 0 \\[3pt]
0 & \frac{1}{2} & \frac{1}{2} \\[3pt]
0 & \frac{1}{2} & \frac{1}{2}
\end{bmatrix}
\times
\begin{bmatrix}
0 & \frac{1}{2} & \frac{1}{2} \\[3pt]
1 & 0 & 0 \\[3pt]
0 & \frac{1}{2} & \frac{1}{2}
\end{bmatrix}.
\end{equation}

We will get back to this observation in Section \ref{Concl}.

%%%%%%%%%%%%%%%%%%%%%%%%%%%%%%%%%%%%%%%%%%%%%%%%%%%%%%%%%%%%%%%
%%%%%%%%%%%%%%%%%%%%%%%%%%%%%%%%%%%%%%%%%%%%%%%%%%%%%%%%%%%%%%%
%%%%%%%%%%%%%%%%%%%%%%%%%%%%%%%%%%%%%%%%%%%%%%%%%%%%%%%%%%%%%%%
%%%%%%%%%%%%%%%%%%%%%%%%%%%%%%%%%%%%%%%%%%%%%%%%%%%%%%%%%%%%%%%

\section{The main result} \label{S:Main-Result}
Our main result, which settles Erd\H{o}s' question for $3 \times 3$ doubly stochastic matrices is the following theorem. The proof is long and processed in the remaining part of the paper.

\begin{thm}\label{thm - maximal}
Let $A$ be a $3\times 3$ doubly stochastic matrix. Then
\[
\|A\|_{\operatorname{F}}^2 \,=\, \max\limits_{\tr}(A)
\]
if and only if there are permutation matrices $P$ and $Q$ such that $PAQ$ is equal to any of the following matrices:
\medskip
\begin{enumerate}[(i)]
\hspace{-5pt}
\begin{minipage}{0.46\linewidth}
\item
\[
\hspace{-5pt}
I_3 \,=\,
		\begin{bmatrix}
			\makebox[\widthof{$\frac{1}{3}$}][l]{$\hspace{.5pt}1$} & 0 & 0 \\[3pt]
			0 & \makebox[\widthof{$\frac{1}{3}$}][l]{$\hspace{.5pt}1$} & 0 \\[3pt]
			0 & 0 & \makebox[\widthof{$\frac{1}{3}$}][l]{$\hspace{.5pt}1$}
		\end{bmatrix};
		\]
		\item
		\[\hspace{-5pt}
		J_3 \,=\,
		\begin{bmatrix}
			\frac{1}{3} &  \frac{1}{3} &  \frac{1}{3} \\[3pt]
			\frac{1}{3} &  \frac{1}{3} &  \frac{1}{3} \\[3pt]
			\frac{1}{3} &  \frac{1}{3} & \frac{1}{3}
		\end{bmatrix};
		\]
		\item
		\[\hspace{-26pt}
		I_1 \oplus J_2 \,=\,
		\begin{bmatrix}
			\makebox[\widthof{$\frac{1}{3}$}][l]{$\hspace{.5pt}1$} & 0 & 0 \\[3pt]
			0 & \tfrac{1}{2} &  \tfrac{1}{2} \\[3pt]
			0 & \tfrac{1}{2} & \tfrac{1}{2}
		\end{bmatrix};
		\]
\end{minipage}
\begin{minipage}{0.08\linewidth}
	\hfill
\end{minipage}
\begin{minipage}{0.46\linewidth}
\item
\[\hspace{-15pt}
S \,=\,
\begin{bmatrix}
	0 & \tfrac{1}{2} & \tfrac{1}{2} \\[3pt]
	\tfrac{1}{2} & \frac{1}{4} & \frac{1}{4} \\[3pt]
	\tfrac{1}{2} & \tfrac{1}{4} & \tfrac{1}{4}
\end{bmatrix};
\]
\item
\[\hspace{-15pt}
T \,=\,
\begin{bmatrix}
	0 & \tfrac{1}{2} & \tfrac{1}{2} \\[3pt]
	\tfrac{1}{2} & 0 & \tfrac{1}{2} \\[3pt]
	\tfrac{1}{2} & \tfrac{1}{2} & 0
\end{bmatrix};
\]
\item
\[\hspace{-15pt}
R \,=\,
\begin{bmatrix}
	\tfrac{3}{5} & 0 & \tfrac{2}{5} \\[3pt]
	0 & \tfrac{3}{5} & \tfrac{2}{5} \\[3pt]
	\tfrac{2}{5} & \tfrac{2}{5} & \tfrac{1}{5}
\end{bmatrix}.
\]
\end{minipage}
\end{enumerate}
\end{thm}

The first three matrices are implied by propositions \ref{ErdosMarcusReePositive} and \ref{ErdosMarcusReeThm}. As for the last two matrices, they are new to this study.

%%%%%%%%%%%%%%%%%%%%%%%%%%%%%%%%%%%%%%%%%%%%%%%%%%%%%%%%%%%%%%%
%%%%%%%%%%%%%%%%%%%%%%%%%%%%%%%%%%%%%%%%%%%%%%%%%%%%%%%%%%%%%%%
%%%%%%%%%%%%%%%%%%%%%%%%%%%%%%%%%%%%%%%%%%%%%%%%%%%%%%%%%%%%%%%
%%%%%%%%%%%%%%%%%%%%%%%%%%%%%%%%%%%%%%%%%%%%%%%%%%%%%%%%%%%%%%%

\section{A weak form of Erd\H{o}s' question}\label{Section3}

\begin{figure}[h]
	\centering
	\begin{subfigure}{.5\textwidth}
		\centering
		\includegraphics[width=.985\textwidth]{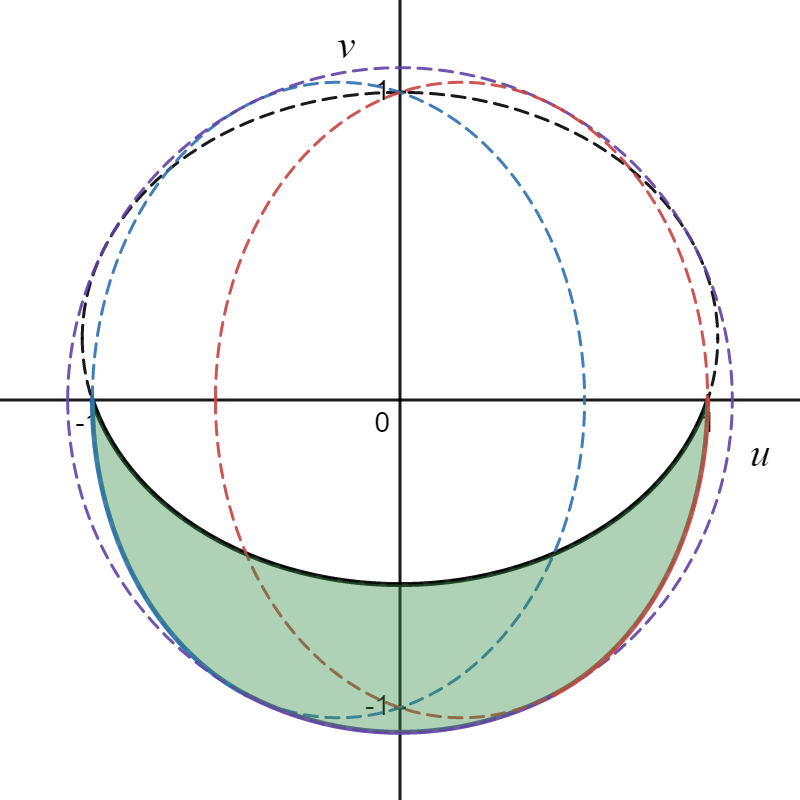}
		\subcaption{The region $\mathcal{U}_-$.}
	\end{subfigure}%
	\begin{subfigure}{.5\textwidth}
		\centering
		\includegraphics[width=.985\textwidth]{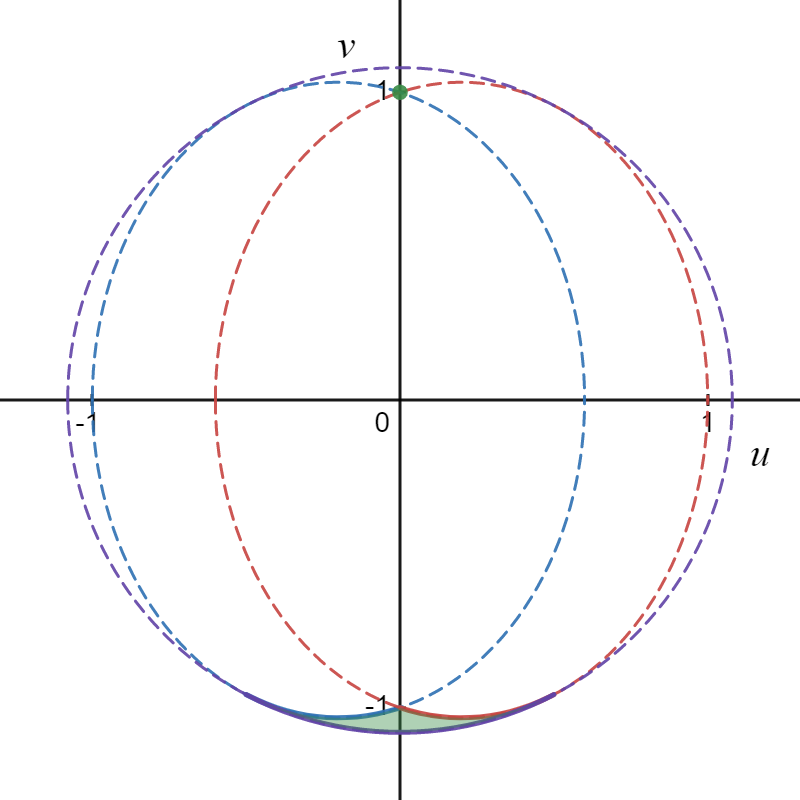}
		\subcaption{The region $\mathcal{U}_+$.}
	\end{subfigure}
	\caption{The region $\mathcal{U}_{-}$ and $\mathcal{U}_+$.}
	\label{figure 1}
\end{figure}

In this section, we characterize the $3 \times 3$ doubly stochastic matrices $A$ such that
\[
\|A\|_{\text{F}}^2 \,=\, \sum\limits_{i=1}^n a_{i, \sigma(i)}
\]
for some permutation $\sigma$. Note that the quantity on the right side of identity is not necessarily $\max\limits_{\tr}(A)$. To do so, we need the following regions: the closed disc
\[
\begin{array}{ll}
\mathcal{E}_{0}: & u^2+v^2 \leq 7/6,
\end{array}
\]
and the following three solid ellipses (boundaries and interiors),
\[
\begin{array}{ll}
\mathcal{E}_{1}: & \frac{(u+1/5)^2}{\, 16/25 \,} + \frac{v^2}{\, 16/15 \,} \leq 1,\\
\mathcal{E}_{2}: & \frac{(u-1/5)^2}{\, 16/25 \,} + \frac{v^2}{\, 16/15 \,} \leq 1, \\
\mathcal{E}_{3}: & \frac{u^2}{\, 16/15 \,} + \frac{(v-1/5)^2}{\, 16/25 \,} \leq 1.
\end{array}
\]
Then we define $\mathcal{U}_{-}$ and $\mathcal{U}_{+}$ as the shaded regions shown in \Cref{figure 1}. In this figure, the color code for the boundary of figures is as follows:
\[
\begin{array}{llll}
\partial \mathcal{E}_{0} \mbox{ in purple}, & \partial \mathcal{E}_{1} \mbox{ in blue}, & \partial \mathcal{E}_{2} \mbox{ in red},  & \partial \mathcal{E}_{3} \mbox{ in black}.
\end{array}
\]

More explicitly, considering the functions
\begin{eqnarray*}
f(u) &\!\!=\!\!& -\sqrt{\frac{3+2|u|-5u^{2}}{3}}, \\
g(u) &\!\!=\!\!& \frac{1-\sqrt{16-15u^{2}}}{5}, \\
h(u) &\!\!=\!\!& \begin{cases}
	-\sqrt{\frac{7}{6}-u^2} \quad &\mathrm{if }  \,|u| \leq \tfrac{1}{2},\\
	\,f(u) \quad &\mathrm{if } \,\tfrac{1}{2}\leq |u| \leq 1,
\end{cases}% f(u) - \mathbf{1}_{\left[-1,1\right]}(2u) \left( f(u) + \sqrt{\frac{7-6u^2}{6}} \right),
\end{eqnarray*}
the regions $\mathcal{U}_-$ and $\mathcal{U}_+$ are defined as follows.
\begin{enumerate}[(i)]\label{def}
\item $\mathcal{U}_-$ is the locus of points $(u,v)$ such that
\[
h(u) \leq v \leq g(u)
\quad \mbox{and} \quad
|u|\leq 1.
\]
\item $\mathcal{U}_+$ is the region of points satisfying
\[
f(u)\leq v \leq g(u)
\quad \mbox{and} \quad
|u|\leq \tfrac{1}{2},
\]
to which we also add the single (isolated) point $(u,v)=(0,1)$.
\end{enumerate}

We are now ready to state the answer to the weak form of Erd\H{o}s' question.

\begin{lem}\label{thm - 1}
Let $A$ be a $3\times 3$ doubly stochastic matrix. Then there is a permutation $\sigma$ such that
\[
\|A\|_{\operatorname{F}}^2 \,=\, \sum\limits_{i=1}^{3} a_{i, \sigma(i)}
\]
if and only if one of the following holds.
\begin{enumerate}[(i)]
\item $A=J_3$.
\item There are permutation matrices $P$ and $Q$ and parameters $(u,v)\in \mathcal{U}_-$ such that
\begin{align*}
PAQ \,=\, \begin{bmatrix}
\tfrac{v+u+3}{4}&\tfrac{1-2v - \sqrt{7-6u^{2}-6v^{2}}}{8}&*\\
&&\\
0&\tfrac{v-u+3}{4}&*\\
&&\\
*&*&*
\end{bmatrix}.
\end{align*}
\item There are permutation matrices $P$ and $Q$ and parameters $(u,v)\in \mathcal{U}_+$ such that
\begin{align*}
PAQ ~=~ \begin{bmatrix}
\tfrac{v+u+3}{4}&\tfrac{1-2v + \sqrt{7-6u^{2}-6v^{2}}}{8}&*\\
&&\\
0&\tfrac{v-u+3}{4}&*\\
&&\\
*&*&*
\end{bmatrix}.
\end{align*}
\end{enumerate}
In $(ii)$ and $(iii)$, under the proposed conditions, the components $*$ are well-defined and uniquely determined such that $PAQ$ becomes a doubly stochastic matrix.
\end{lem}

The rest of this section is devoted to providing a complete and detailed proof of this result. We start with some reductions.

First, since for any matrix $A$ and for any permutation matrix $P$, we have $\|A\|_{\text{F}}=\|AP\|_{\text{F}}$, instead of the most general identity  $\|A\|_{\text{F}}^2=\tr(AP)$, we may just concentrate on
\begin{equation}\label{E:reduced-weak-erdos}
\|A\|_{\text{F}}^2=\tr(A).
\end{equation}
Second, according to Proposition \ref{ErdosMarcusReePositive}, except for the notable matrix $J_3$, there is no $3\times 3$ doubly stochastic all of whose entries are strictly positive which saturates the inequality \eqref{E:EMR}. Hence, we shall focus our attention to those doubly stochastic matrices $A$ with at least one zero element, and, again thanks to the permutation invariance of the equation, without any loss of generality we assume that $a_{2,1}= 0$. Therefore, the matrix $A$ is of the form
\begin{equation}\label{E:def-A1}
A \,=\,
\begin{bmatrix}
a & c & 1-a-c \\
&&\\
0 & b & 1-b \\
&&\\
1-a & 1-b-c & a+b+c-1
\end{bmatrix},
\end{equation}
where the parameters $a$, $b$ and $c$ are necessarily restricted to a region in $\mathbb{R}^3$ so that the matrix $A$ becomes doubly stochastic. We can continue our calculation with this format. However, the region we obtain is at the intersection of discs and ellipses whose main axes are not parallel to the coordinate axes, and whose center is not the origin. To overcome this feature in advance, we consider the following format
\begin{equation}\label{E:def-A2}
A \,=\,
\begin{bmatrix}
\tfrac{v+u+3}{4} & w & \tfrac{1-v-u}{4}-w \\
&&\\
0 & \tfrac{v-u+3}{4} & \tfrac{1-v+u}{4} \\
&&\\
\tfrac{1-v-u}{4} & \tfrac{1-v+u}{4}-w & \tfrac{v+1}{2}+w
\end{bmatrix},
\end{equation}
which is as general as the previous case (just linear transformation of the parameters), and of course with a different set of conditions on $u$, $v$ and $w$ to ensure that we end up with a doubly stochastic matrix.  But, even though it is not transparent at this stage, we will see that the conditions on the new set of parameters are more symmetric and simpler.

Under the format \eqref{E:def-A2}, after some elementary calculations, the equation \eqref{E:reduced-weak-erdos} simplifies to
\[
4w^{2}+\left(2v-1\right)w+\frac{3u^{2}+5v^{2}-2v-3}{8} \,=\, 0.
\]
Hence, $w$ is given more explicitly by
\begin{equation}\label{E:w-plus-minus}
w \,=\, \frac{1-2v \pm \sqrt{7-6u^{2}-6v^{2}}}{8}.
\end{equation}
Since $w$ ought to be a real number, we need to have
\[
7-6u^{2}-6v^{2} \,\geq\, 0,
\]
i.e., $(u,v)\in \mathcal{E}_{0}$, the closed disc of radius $\smash{\sqrt{7/6}}$ centered at the origin.

We proceed to complete the proof of Lemma \ref{thm - 1}. To this end, we need to identify all the constraints that the parameters $u$ and $v$ and $w$ (or equivalently, the parameters $a$ and $b$ and $c$) must verify so that the matrix $A$ to be doubly stochastic. At first glance on \eqref{E:def-A1}, we see that a necessary and sufficient condition for $A$ to be doubly stochastic is that
\[
\begin{array}{ccc}
0 \leq a \leq 1, &\quad 0 \leq c \leq 1,  \quad& 0 \leq 1-a-c \leq 1,  \\
&&\\
 &\quad 0 \leq b \leq 1,  \quad& 0 \leq 1-b \leq 1,  \\
&&\\
0 \leq 1-a \leq 1,  &\quad 0 \leq 1-b-c \leq 1,  \quad& 0 \leq a+b+c-1 \leq 1.
\end{array}
\]
However, the above set of restrictions simplifies to
\[
\begin{array}{ccc}
 a \,\geq\, 0, &\quad b \,\geq\, 0,  \quad& c \,\geq\, 0,  \\
&&\\
a+c \,\leq\, 1,  &\quad b+c \,\leq\, 1,  \quad&  a+b+c \,\geq\, 1.
\end{array}
\]
By \eqref{E:w-plus-minus}, as there are two possibilities for the values $w$, we also consider two cases. Recall that we already imposed
\begin{equation}\label{E:first-cond}
(u,v)\in \mathcal{E}_{0}.
\end{equation}

\noindent {\underline{\bf Case 1:}} Let
\[
w \,=\, \frac{1-2v - \sqrt{7-6u^{2}-6v^{2}}}{8}, \qquad (u,v)\in \mathcal{E}_{0}.
\]
In order to be a doubly stochastic matrix, there are several restrictions on the elements of $A$ which we treat one by one. In the following we address $A$ according to the recipe \eqref{E:def-A1} or \eqref{E:def-A2}.

\noindent {\bf Condition $a \geq 0$:} Since
\[
a \,=\, \frac{v+u+3}{4},
\]
the condition $a \geq 0$ holds if and only if
\[
u+v \,\geq\, -3.
\]

\noindent {\bf Condition $b \geq 0$:} Since
\[
b \,=\, \frac{v-u+3}{4},
\]
the condition $b \geq 0$ holds if and only if
\[
u-v \,\leq\, 3.
\]

\noindent Remark: Since $\mathcal{E}_{0}$ is a closed disc of radius $\sqrt{7/6}$ and, according to \eqref{E:first-cond}, we require $(u,v)$ be in this closed disc, then
\[
u+v \,\geq\, -2\sqrt{7/6} \,>\, -3,
\]
and
\[
u-v \,\leq\, 2 \sqrt{7/6} \,<\, 3.
\]
In other words, in the light of the constraint $(u,v)\in \mathcal{E}_{0}$, the above two conditions are automatically fulfilled.

\noindent {\bf Condition $c \geq 0$:} Since
\[
c=w \,=\, \frac{1-2v - \sqrt{7-6u^{2}-6v^{2}}}{8},
\]
the required condition $c \geq 0$  is realized if and only if
\[
1-2v \,\geq\, \sqrt{7-6u^{2}-6v^{2}}.
\]
Hence, $v \leq 1/2$ and $(u,v)\in \mathcal{E}_{0}$, and then, if we square both sides and simplify, we obtain
\[
\frac{u^2}{\, 16/15 \,} + \frac{(v-1/5)^2}{\, 16/25 \,} \,\geq\, 1,
\]
which means $(u,v) \notin \Int(\mathcal{E}_3)$.

\noindent {\bf Condition $a+c \leq 1$:} Since
\[
a+c \,=\, \frac{v+u+3}{4}+w \,=\, \frac{7+2u-\sqrt{7-6u^{2}-6v^{2}}}{8},
\]
the condition $a+c \leq 1$ holds if and only if
\[
1-2u+\sqrt{7-6u^{2}-6v^{2}} \,\geq\, 0.
\]
This happens if and only if
\begin{enumerate}[-]
\item either $u \leq 1/2$ and $(u,v)\in \mathcal{E}_{0}$,
\item or $u \geq 1/2$ and
\[
\sqrt{7-6u^{2}-6v^{2}} \,\geq\, 2u-1.
\]
After squaring and simplifying, we obtain
\[
\frac{(u-1/5)^2}{\, 16/25 \,} + \frac{v^2}{\, 16/15 \,} \,\leq\, 1,
\]
which means $(u,v)\in \mathcal{E}_2$.
\end{enumerate}

\noindent {\bf Condition $b+c \leq 1$:} Since
\[
b+c \,=\, \frac{v-u+3}{4}+w  \,=\, \frac{7-2u-\sqrt{7-6u^{2}-6v^{2}}}{8}
\]
the condition $b+c \leq 1$ holds if and only if
\[
1+2u+\sqrt{7-6u^{2}-6v^{2}} \,\geq\, 0.
\]
This happens if and only if
\begin{enumerate}[-]
\item either $u \geq -1/2$ and $(u,v)\in \mathcal{E}_{0}$,
\item or $u \leq -1/2$ and
\[
\sqrt{7-6u^{2}-6v^{2}} \,\geq\, -1-2u.
\]
After squaring and simplifying, we obtain
\[
\frac{(u+1/5)^2}{\, 16/25 \,} + \frac{v^2}{\, 16/15 \,} \,\leq\, 1,
\]
which means $(u,v)\in \mathcal{E}_{1}$
\end{enumerate}

\noindent {\bf Condition $a+b+c\geq 1$:} Since
\[
a+b+c \,=\, \tfrac{v+u+3}{4}+ \tfrac{v-u+3}{4}+w \,=\, \frac{2v+13 - \sqrt{7-6u^{2}-6v^{2}}}{8}
\]
the condition $a+b+c \geq 1$ holds if and only if
\[
2v+5-\sqrt{7-6u^{2}-6v^{2}} \,\geq\, 0.
\]
But, since $(u,v)\in \mathcal{E}_{0}$, we have
\[
2v+5-\sqrt{7-6u^{2}-6v^{2}} \,\geq\, -2\sqrt{7/6}+5 - \sqrt{7} \approx0.19 \,>\,  0.
\]
This condition is also automatically satisfied.

\begin{figure}[H]
    \centering
    \includegraphics[width=\textwidth]{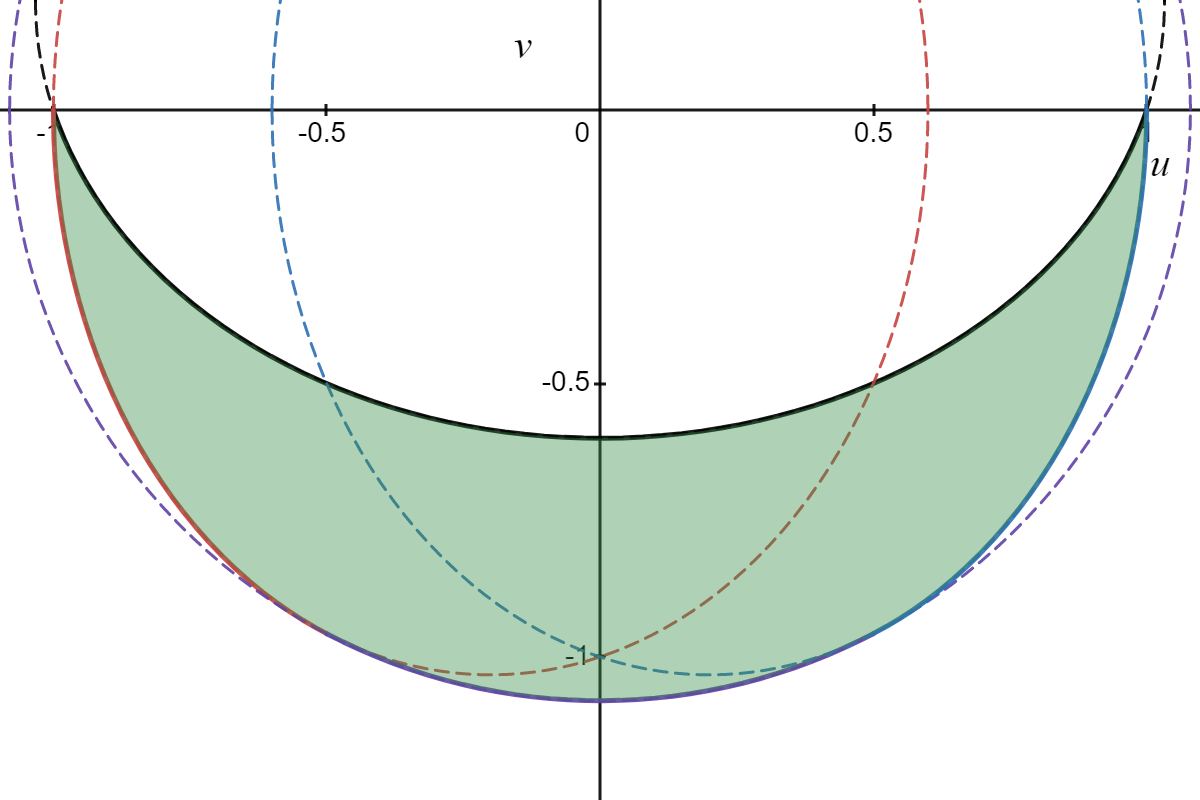}
    \caption{The region $\mathcal{U}_-$.}\label{fig:my_label1}
\end{figure}

\noindent {\bf Summary of the Case 1:}
The matrix
\[
A ~=~
\begin{bmatrix}
\tfrac{v+u+3}{4} & w & \tfrac{1-v-u}{4}-w \\
&&\\
0 & \tfrac{v-u+3}{4} & \tfrac{1-v+u}{4} \\
&&\\
\tfrac{1-v-u}{4} & \tfrac{1-v+u}{4}-w & \tfrac{v+1}{2}+w
\end{bmatrix},
\]
where $w=\frac{1-2v - \sqrt{7-6u^{2}-6v^{2}}}{8}$, is doubly stochastic and satisfies
\[
\sum_{i=1}^{3}\sum_{j=1}^{3} a_{i,j}^2 = \sum_{i=1}^{3} a_{i,i}
\]
if and only if the following conditions simultaneously hold:
\begin{enumerate}[(i)]
\item $(u,v)\in \mathcal{E}_{0}$,
\item $(u,v) \notin \Int(\mathcal{E}_3)$,
\item $v\leq 1/2$,
\item if $u \geq 1/2$, $(u,v)\in \mathcal{E}_2$,
\item if $u\leq -1/2$, $(u,v)\in \mathcal{E}_1$.
\end{enumerate}
More explicitly, the resulting region can also be described as follows: it is the set of points $(u,v)$ enclosed by the graph of functions $h$ and $g$ in the range $|u|\leq 1$, where
\begin{gather*}
    g(u) \,=\, \frac{1-\sqrt{16-15u^{2}}}{5},
    \quad\&\quad
    h(u) \,=\, \begin{cases}
            -\sqrt{\frac{3+2|u|-5u^{2}}{3}} \quad &\mathrm{if } \,\tfrac{1}{2}\leq |u| \leq 1;\\
            -\sqrt{\frac{7-6u^2}{6}} \quad &\mathrm{if }  \,|u| \leq \tfrac{1}{2}.
    \end{cases}
\end{gather*}
See Figure \ref{fig:my_label1}.

\noindent In Figure \ref{fig:my_label1}, the boundary of $\mathcal{E}_{0}$ is in purple, whereas the boundaries of $\mathcal{E}_1$, $\mathcal{E}_2$ and $\mathcal{E}_3$ are respectively colored in blue, red and black.\\

\noindent {\underline{\bf Case 2:}} In this case we consider
\[
w\,=\,\frac{1-2v + \sqrt{7-6u^{2}-6v^{2}}}{8}.
\]
As in the Case 1, in order to have a doubly stochastic matrix, there are several constraints on the elements of $A$. As in the previous case, the conditions $a\geq0$ and $b\geq 0$ are automatically fulfilled. Therefore, we treat the remaining four needed conditions below.

\noindent {\bf Condition $c \geq 0$:} Since
\[
c=w \,=\, \frac{1-2v + \sqrt{7-6u^{2}-6v^{2}}}{8},
\]
the required condition $c \geq 0$  is realized if and only if
\[
1-2v + \sqrt{7-6u^{2}-6v^{2}} \,\geq\, 0.
\]
This is similar to Condition $a+c\leq 1$ in Case 1 with $u$ and $v$ exchanged so this happens if and only if
\begin{enumerate}[-]
\item either $v \leq 1/2$ and $(u,v)\in \mathcal{E}_{0}$,
\item or $v \geq 1/2$ and $(u,v)\in \mathcal{E}_3$.
\end{enumerate}

\noindent {\bf Condition $a+c \leq 1$:} Since
\[
a+c \,=\, \frac{v+u+3}{4}+w \,=\, \frac{7+2u+\sqrt{7-6u^{2}-6v^{2}}}{8},
\]
the condition $a+c \leq 1$ holds if and only if
\[
1-2u \,\geq\, \sqrt{7-6u^{2}-6v^{2}}.
\]
This is analogous to Condition $c\geq 0$ in Case 1 with $u$ and $v$ exchanged, and thus this happens if and only if $u \leq 1/2$, $(u,v)\in \mathcal{E}_{0}$ and $(u,v) \notin \Int(\mathcal{E}_2)$.

\noindent {\bf Condition $b+c \leq 1$:} Since
\[
b+c \,=\, \frac{v-u+3}{4}+w  \,=\, \frac{7-2u+\sqrt{7-6u^{2}-6v^{2}}}{8},
\]
the condition $b+c \leq 1$ holds if and only if
\[
1+2u \,\geq\, \sqrt{7-6u^{2}-6v^{2}}.
\]
Once again, this is similar to Condition $c\geq 0$ in Case 1 with $(u,v)$ replaced by $(v,-u)$ so this happens if and only if $u \geq -1/2$, $(u,v)\in \mathcal{E}_{0}$ and $(u,v) \notin \Int(\mathcal{E}_1)$.

\noindent {\bf Condition $a+b+c\geq 1$:} Since
\[
a+b+c \,=\, \tfrac{v+u+3}{4}+ \tfrac{v-u+3}{4}+w \,=\, \frac{2v+13 + \sqrt{7-6u^{2}-6v^{2}}}{8},
\]
the condition $a+b+c \geq 1$ holds if and only if
\[
2v+5+\sqrt{7-6u^{2}-6v^{2}} \,\geq\, 0.
\]
But, since $(u,v)\in \mathcal{E}_{0}$, we have
\[
2v+5+\sqrt{7-6u^{2}-6v^{2}} \,\geq -2\sqrt{7/6}+5 \approx 2.84 \,>\,  0.
\]
Thus, this condition is also automatically satisfied.

\noindent {\bf Summary of the Case 2:}
The matrix
\[
A \,=\,
\begin{bmatrix}
	\tfrac{v+u+3}{4} & w & \tfrac{1-v-u}{4}-w \\
	&&\\
	0 & \tfrac{v-u+3}{4} & \tfrac{1-v+u}{4} \\
	&&\\
	\tfrac{1-v-u}{4} & \tfrac{1-v+u}{4}-w & \tfrac{v+1}{2}+w
\end{bmatrix},
\]
where $w=\frac{1-2v + \sqrt{7-6u^{2}-6v^{2}}}{8}$, is doubly stochastic and satisfies
\[
\sum_{i=1}^{3}\sum_{j=1}^{3} a_{i,j}^2 \,=\, \sum_{i=1}^{3} a_{i,i}
\]
if and only if the following conditions simultaneously hold:
\begin{enumerate}[(i)]
	\item $(u,v)\in \mathcal{E}_{0}$;
	\item $|u| \leq \tfrac{1}{2}$;
	\item $(u,v)\notin \Int(\mathcal{E}_1) $;
	\item $(u,v)\notin \Int(\mathcal{E}_2)$.
	\item if $v\geq 1/2$, $(u,v)\in\mathcal{E}_3$.
\end{enumerate}
The resulting region can be described as follows: it is the set of points $(u,v)$ enclosed by the graph of functions $f$ and $g$ in the range $|u|\leq \tfrac{1}{2}$, where
\begin{gather*}
    g(u)\,:=~ \frac{1-\sqrt{16-15u^{2}}}{5},
    \quad\&\quad
    f(u)\,:=\, \smash{-\sqrt{\frac{3+2\left|u\right|-5u^{2}}{3}}}.
\end{gather*}
along with the isolated point $(0,1)$. Together they form the following region.

\begin{figure}[H]
    \centering
    \includegraphics[width=\textwidth]{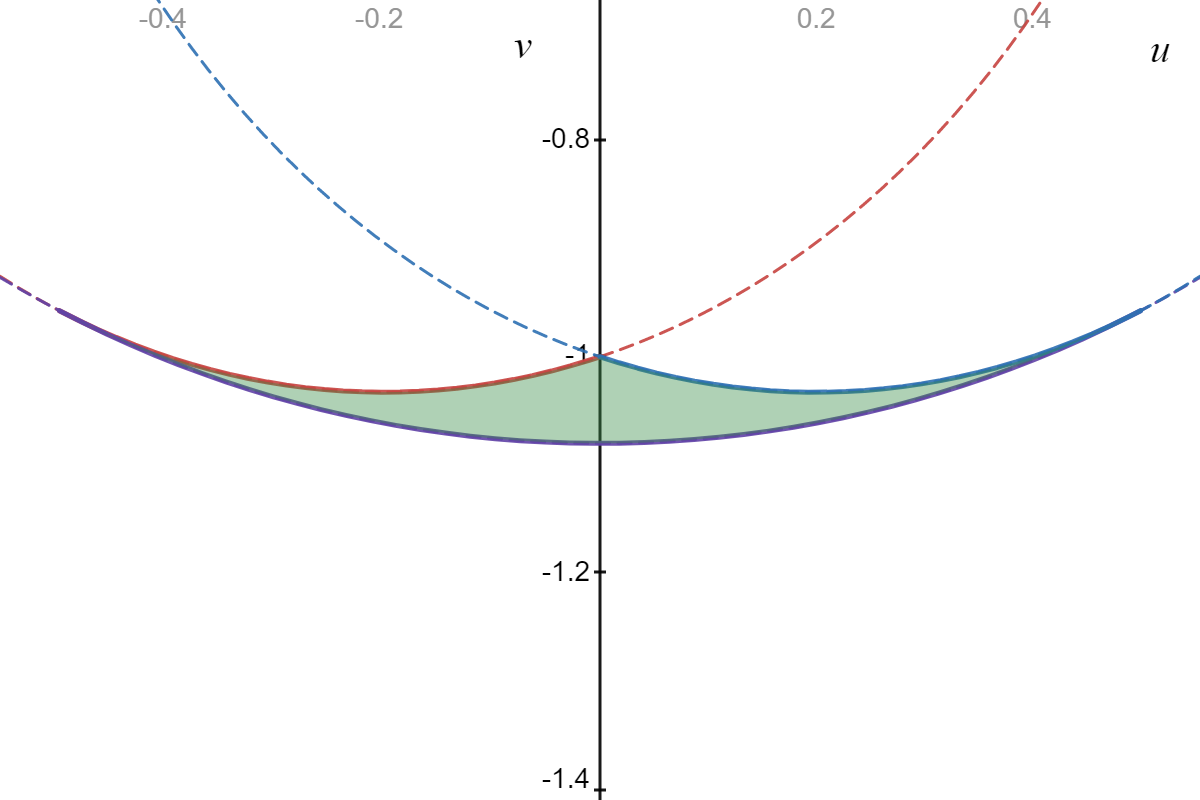}
    \caption{The region $\mathcal{U}_+ \!\setminus\! (0,1)$.}    \label{fig:my_label2}
\end{figure}

\noindent In Figure \ref{fig:my_label2}, the boundary of $\mathcal{E}_{0}$ is in purple, whereas the boundaries of $\mathcal{E}_1$ and $\mathcal{E}_2$ are respectively represented in blue and red.\\

\noindent Remark: With the notable exception of the point $(u,v) = (0,1)$, the set $\mathcal{U}_+$ is a subset of $\mathcal{U}_-$. However, this should not be misleading as the inclusion does not imply that the set of doubly stochastic matrices described in Case 2 is included in the one described in Case 1. In fact, we can show that they are disjoint sets in $\mathbb{R}^9$.

%%%%%%%%%%%%%%%%%%%%%%%%%%%%%%%%%%%%%%%%%%%%%%%%%%%%%%%%%%%%%%%
%%%%%%%%%%%%%%%%%%%%%%%%%%%%%%%%%%%%%%%%%%%%%%%%%%%%%%%%%%%%%%%
%%%%%%%%%%%%%%%%%%%%%%%%%%%%%%%%%%%%%%%%%%%%%%%%%%%%%%%%%%%%%%%
%%%%%%%%%%%%%%%%%%%%%%%%%%%%%%%%%%%%%%%%%%%%%%%%%%%%%%%%%%%%%%%

\section{The proof of Theorem \ref{thm - maximal}} \label{S:Prrof-Main-Result}
Since $J_3$ trivially saturates the inequality \ref{E:EMR}, in the remainder of this section we assume that $A \neq J_3$. Of course, our main assumption is that  $\|A\|_{\text{F}}^2 = \max\limits_{\tr}(A)$. Hence, $A$ is necessarily of the form \eqref{E:def-A2},
where the parameters satisfy the restrictions described in Lemma \ref{thm - 1}. However, not all such matrices saturate the inequality \ref{E:EMR}. For example, choosing $u = 0$ and $v = -\tfrac{21}{20}$ give us a doubly stochastic matrix $A$ such that
\[
\|A\|_{\text{F}}^2 \,=\, \tr(A) \,<\, \max\limits_{\tr}(A).
\]
In short, the $3 \times 3$ doubly stochastic matrices that are solutions to the equation $\|A\|_{\text{F}}^2 = \max\limits_{\tr}(A)$ form a proper subset of the regions characterized in Lemma \ref{thm - 1}. In fact, this subset corresponds to those matrices that also verify the additional condition $\tr(AP) \leq \tr(A)$ for all permutations $P$. As there are only 5 non-trivial permutation matrices of order $3 \times 3$, namely
\begin{gather*}
    P_1 \,=\,
    \begin{bmatrix}
       0 & 0 & 1 \\
       1 & 0 & 0 \\
       0 & 1 & 0
    \end{bmatrix},
    \quad
    P_2 \,=\,
      \begin{bmatrix}
       0 & 1 & 0 \\
       0 & 0 & 1 \\
       1 & 0 & 0
    \end{bmatrix},
    \quad
    P_3 \,=\,
    \begin{bmatrix}
       0 & 1 & 0 \\
       1 & 0 & 0 \\
       0 & 0 & 1
    \end{bmatrix},\\
    P_4 \,=\,
    \begin{bmatrix}
       0 & 0 & 1 \\
       0 & 1 & 0 \\
       1 & 0 & 0
    \end{bmatrix},
    \quad
    P_5 \,=\,
    \begin{bmatrix}
       1 & 0 & 0 \\
       0 & 0 & 1 \\
       0 & 1 & 0
    \end{bmatrix},
\end{gather*}
this means that we need to check the condition $\tr(AP) \leq \tr(A)$ for the above five permutations and see if they give rise to new restrictions.

\noindent {\bf Permutation $P_1$:} We have $\tr(A) \geq \tr(AP_1)$
\begin{eqnarray}
&\iff& (a)+(b)+(a+b+c-1) \geq (c)+(1-b)+(1-a) \notag\\
&\iff& a+b \geq 1 \notag\\
&\iff& \frac{v+u+3}{4}+\frac{v-u+3}{4} \geq1 \notag\\
&\iff& v \geq -1. \label{NewConstraint1j}
\end{eqnarray}

\noindent {\bf Permutation $P_2$:} We have $\tr(A) \geq \tr(AP_2)$
\begin{eqnarray}
	&\iff& (a)+(b)+(a+b+c-1) \geq (1-a-c)+(0)+(1-b-c) \notag\\
	&\iff& a+b+c \geq 1. \label{NewConstraint2j}
\end{eqnarray}

\noindent {\bf Permutation $P_3$:} We have $\tr(A) \geq \tr(AP_3)$
\begin{eqnarray}
&\iff& (a)+(b)+(a+b+c-1) \geq (c)+(0)+(a+b+c-1) \notag\\
&\iff& a+b \geq c.\label{NewConstraint3j}
\end{eqnarray}

\noindent {\bf Permutation $P_4$:} We have $\tr(A) \geq \tr(AP_4)$
\begin{eqnarray}
&\iff& (a)+(b)+(a+b+c-1) \geq (1-a-c)+(b)+(1-a) \notag\\
&\iff& 4a+b+2c \geq 3 \notag\\
&\iff& (v+u+3)+\frac{v-u+3}{4} + \frac{1-2v \pm\sqrt{7-6u^{2}-6v^{2}}}{4} \geq 3 \notag\\
&\iff& 4+3v+3u \pm \sqrt{7-6u^{2}-6v^{2}} \geq 0. \label{NewConstraint4j}
\end{eqnarray}

\noindent {\bf Permutation $P_5$:} We have $\tr(A) \geq \tr(AP_5)$
\begin{eqnarray}
&\iff& (a)+(b)+(a+b+c-1) \geq (a)+(1-b)+(1-b-c) \notag\\
&\iff& a+4b+2c \geq 3 \notag\\
&\iff& \frac{v+u+3}{4}+(v-u+3) + \frac{1-2v \pm\sqrt{7-6u^{2}-6v^{2}}}{4} \geq 3 \notag\\
&\iff& 4+3v-3u \pm \sqrt{7-6u^{2}-6v^{2}} \geq 0. \label{NewConstraint5j}
\end{eqnarray}
We can summarize \eqref{NewConstraint4j} and \eqref{NewConstraint5j} as
\begin{equation}\label{NewConstraint6j}
4+3v-3|u| \pm \sqrt{7-6u^{2}-6v^{2}} \geq 0.
\end{equation}

Now, observe that \eqref{NewConstraint2j} is automatically satisfied since $a+b+c\geq1$ is one of the conditions that was respected in the proof of Lemma \ref{thm - 1}. Assuming \eqref{NewConstraint1j}, then \eqref{NewConstraint3j} is also automatically satisfied since then
\[
a+b \,\geq\, 1 \,\geq\, c.
\]
Therefore, we only need to asses the restrictions imposed by the conditions \eqref{NewConstraint1j} and \eqref{NewConstraint6j}. To do so, we consider in the following two cases corresponding respectively to the $\pm$ cases.

\noindent {\bf Case 1, The Minus Sign:} The requirement \eqref{NewConstraint6j} becomes
\[
4+3v-3|u| \,\ge\, \sqrt{7-6u^{2}-6v^{2}},
\]
which is fulfilled if and only if $(u,v)\in \mathcal{E}_{0}$ (which is trivially satisfied) and
\begin{equation}\label{eq - important}
\begin{cases}  (4+3v-3|u|)^2 +6u^{2}+6v^{2}-7 \,\geq\, 0, \\ 4+3v-3|u| \,\geq\, 0. \end{cases}
\end{equation}
Since the function in the first inequality is always convex relative to $v$ for any given $u$, its maximum on a given interval $v\in[\alpha,\beta]$ is always attained at one of the endpoints. 

According to Lemma \ref{thm - 1}, in this case $(u,v)\in \mathcal{U}_-$, which implies $v\leq \frac{1-\sqrt{16-15u^{2}}}{5}$. Hence, respecting \eqref{NewConstraint1j}, we should have
\[
-1 \leq v \leq \frac{1-\sqrt{16-15u^{2}}}{5}.
\]
This condition can be incorporated with the second inequality of \eqref{eq - important} to give
\begin{equation}\label{NewConstraint}
	\max\left(-1,|u|-\tfrac{4}{3}\right) \leq v \leq \frac{1-\sqrt{16-15u^{2}}}{5}.
\end{equation}
Therefore, the function $(4+3v-3|u|)^2 +6u^{2}+6v^{2}-7$ is always bounded below by the maximum of its evaluation at the two extremal values of $v$ in \eqref{NewConstraint}.

Evaluating the lower bound $\max\left(-1,|u|-\tfrac{4}{3}\right)$ gives
\[
\begin{cases}
\left(5\left|u\right|-2\right)\left|u\right|, \qquad &\text{if } |u|\leq \tfrac{1}{3}; \\
4u^{2}-\frac{16}{3}\left|u\right|+\frac{11}{9}, \qquad &\text{if } |u|\geq \tfrac{1}{3},
\end{cases}
\]
which is always strictly negative on $\mathcal{U}_-$, with the exception of $(u,v)=(0,-1)$.

Evaluating the upper bound $\frac{1-\sqrt{16-15u^{2}}}{5}$ gives

\[
g(u) \,:=~ \frac{2}{5}\left(5u^{2}-23\left|u\right|+20-\left(5-3\left|u\right|\right)\sqrt{16-15u^{2}}\right).
\]
Since $5u^{2}-23\left|u\right|+20$ is positive for $(u,v)\in\mathcal{U}_-$, $g(u)$ is positive if and only if
\begin{align}\label{eq - important2}
\left(20-23\left|u\right|+5u^{2}\right)^{\!2}\!-\left(5-3\left|u\right|\right)^{2}\left(16-15u^{2}\right) ~\geq~ 0,
\end{align}
which simplifies to
\[
40\left|u\right|\left(\left|u\right|-1\right)\left(4u^{2}-13\left|u\right|+11\right) ~\ge~ 0.
\]
However, since $4u^{2}-13\left|u\right|+11$ is always positive on $\mathcal{U}_-$ and since $|u|\leq1$ for every $(u,v)\in\mathcal{U}_-$, the above inequality is never verified, except at $u=0$ and $u=\pm 1$. Consequently, \cref{eq - important2} is never fulfilled and $g(u)\leq0$, except when $u=0$ and $u=\pm 1$.
Therefore, we have
\[
(4+3v-3|u|)^2 +6u^{2}+6v^{2}-7 \,\leq\, \max\{f(u),\,g(u)\} ~\leq~ 0,
\]
with equality only if $(u,v)$ is equal to one of $(0,-1),(0,-3/5),(1,0)$ or $(-1,0)$. It then follows by a direct verification that these points are indeed solutions to the constraints \eqref{eq - important}, and thus that they correspond to the following  doubly stochastic matrices:
\[
\begin{bmatrix}
   \tfrac{1}{2}&\tfrac{1}{4}&\tfrac{1}{4}\\[3pt]
   0&\tfrac{1}{2}&\tfrac{1}{2}\\[3pt]
   \tfrac{1}{2}&\tfrac{1}{4}&\tfrac{1}{4}
\end{bmatrix},
\quad
\begin{bmatrix}
   \tfrac{3}{5}&0&\tfrac{2}{5}\\[3pt]
   0&\tfrac{3}{5}&\tfrac{2}{5}\\[3pt]
   \tfrac{2}{5}&\tfrac{2}{5}&\tfrac{1}{5}
\end{bmatrix},
\quad
\begin{bmatrix}
   1&0&0\\[3pt]
   0&\tfrac{1}{2}&\tfrac{1}{2}\\[3pt]
   0&\tfrac{1}{2}&\tfrac{1}{2}
\end{bmatrix}
\quad\&\quad
\begin{bmatrix}
   \tfrac{1}{2}&0&\tfrac{1}{2}\\[3pt]
   0&1&0\\[3pt]
   \tfrac{1}{2}&0&\tfrac{1}{2}
\end{bmatrix}.
\]
The first matrix is permutationally equal to $S$. The second is $R$, the third is $I_1 \oplus J_2$, and the last one is permutationally equal to $I_1 \oplus J_2$.\\

\noindent {\bf Case 2, The Plus Sign:} The requirement \eqref{NewConstraint6j} becomes
\begin{align}\label{eq - extra}
4+3v-3|u| + \sqrt{7-6u^{2}-6v^{2}} ~\ge~ 0.
\end{align}
In the proof of \Cref{thm - 1}, we saw that whenever $v\geq 0$, the only pair $(u,v)$ in $\mathcal{U}_+$ is the isolated pair $(0,1)$. It is then easily verified that
\[
8 ~=~ 4+3(1)-3|0| + \sqrt{7-6(0)^{2}-6(1)^{2}} ~\ge~ 0,
\]
and thus that $(u,v)=(0,1)$ will indeed provide a doubly stochastic solution to the equation $\|A\|_{\text{F}}^2=\max\limits_{\tr}(A)$.

Let us now consider the case $v\leq 0$ and observe that \eqref{NewConstraint1j} and \Cref{thm - 1} ensures that
\begin{align}\label{eq - 3}
-1 \leq v \leq -\sqrt{\frac{3+2\left|u\right|-5u^{2}}{3}}.
\end{align}
Consequently, we have $1 \geq \frac{3+2\left|u\right|-5u^{2}}{3}$ which simplifies to $0\ge\left|u\right|\left(2-5\left|u\right|\right)$. Hence, we either have $u=0$ or $|u|\geq 2/5$. If $u=0$, then \eqref{eq - 3} implies that $v=-1$, and we directly verify that this is indeed a doubly stochastic solution to the equation $\|A\|_{\text{F}}^2=\max\limits_{\tr}(A)$.

Otherwise, suppose that $|u|\geq 2/5$ and recall that we also have $|u|\leq 1/2$ because of \Cref{thm - maximal}. In this case, observe that
\[
\frac{\partial }{\partial v} \left( 4+3v-3|u| + \sqrt{7-6u^{2}-6v^{2}} \right) \,=\, 3-\frac{6v}{\sqrt{7-6u^{2}-6v^{2}}} ~\geq~ 0
\]
since $v\leq 0$. It then follows from \eqref{eq - 3} that
\begin{align*}
    0 ~&\leq~ 4+3v-3|u| + \sqrt{7-6u^{2}-6v^{2}} \\
    &\leq~ 4-3\sqrt{\frac{3+2\left|u\right|-5u^{2}}{3}}-3|u| + \sqrt{7-6u^{2}-6\left( -\sqrt{\tfrac{3+2\left|u\right|-5u^{2}}{3}} \right)^{\!\!2}} \\
    &=~ 5-5\left|u\right|-\sqrt{3\left(3+2\left|u\right|-5u^{2}\right)}.
\end{align*}
However, since $2/5 \leq |u| \leq 1/2$, $5-5\left|u\right|$ is always positive and thus
\begin{align*}
    0 \leq 5-5\left|u\right|-\sqrt{3\left(3+2\left|u\right|-5u^{2}\right)} ~&\iff~
    \left(5-5\left|u\right|\right)^{2} \ge 3\left(3+2\left|u\right|-5u^{2}\right) \\
    &\iff~ \left(2-5\left|u\right|\right)\left(1-\left|u\right|\right)\ge0.
\end{align*}
Consequently, since we supposed that $2/5 \leq |u| \leq 1/2$, the only possible solutions occurs at $|u|=2/5$. A direct verification then suffice to ensure that these are indeed doubly stochastic solutions to the equation $\|A\|_{\text{F}}^2=\max\limits_{\tr}(A)$

Evaluating at $(u,v)=(0,1), (0,-1), (-2/5,-1)$, and $(2/5,-1)$ respectively give us the following matrices:
\begin{gather*}
    \begin{bmatrix}
       1&0&0\\[3pt]
       0&1&0\\[3pt]
       0&0&1
    \end{bmatrix},
    \quad
    \begin{bmatrix}
       \tfrac{1}{2}&\tfrac{1}{2}&0\\[3pt]
       0&\tfrac{1}{2}&\tfrac{1}{2}\\[3pt]
       \tfrac{1}{2}&0&\tfrac{1}{2}
    \end{bmatrix},
    \quad
    \begin{bmatrix}
       \tfrac{2}{5}&\tfrac{2}{5}&\tfrac{1}{5}\\[3pt]
       0&\tfrac{3}{5}&\tfrac{2}{5}\\[3pt]
       \tfrac{3}{5}&0&\tfrac{2}{5}
    \end{bmatrix}
    \quad\&\quad
    \begin{bmatrix}
       \tfrac{3}{5}&\tfrac{2}{5}&0\\[3pt]
       0&\tfrac{2}{5}&\tfrac{3}{5}\\[3pt]
       \tfrac{2}{5}&\tfrac{1}{5}&\tfrac{2}{5}
    \end{bmatrix}.
\end{gather*}
The first matrix is $I_3$, the second is permutationally equal to $T$. The last two are permutationally equal to $R$. We have thus identified a certain number of matrices verifying the set of constraints ensuring that the inequality \eqref{E:EMR} is saturated. To complete the demonstration of \Cref{thm - maximal}, one only need to notice that both the Frobenius norm and the maximal trace are permutation invariant functions.

%%%%%%%%%%%%%%%%%%%%%%%%%%%%%%%%%%%%%%%%%%%%%%%%%%%%%%%%%%%%%%%
%%%%%%%%%%%%%%%%%%%%%%%%%%%%%%%%%%%%%%%%%%%%%%%%%%%%%%%%%%%%%%%
%%%%%%%%%%%%%%%%%%%%%%%%%%%%%%%%%%%%%%%%%%%%%%%%%%%%%%%%%%%%%%%
%%%%%%%%%%%%%%%%%%%%%%%%%%%%%%%%%%%%%%%%%%%%%%%%%%%%%%%%%%%%%%%

\section{Concluding remarks}\label{Concl}

\begin{enumerate}
\item We thank the anonymous reviewers for their careful reading of our manuscript and their many insightful comments. In particular, one of the reviewer pointed out to us that the matrix $T$ in Theorem 4.1 can be generalized to any order $n \geq 2$. Indeed, let
\[
T_n :=~ \frac{1}{n-1} (n J_n-I_n).
\]
Then
\[
\|T_n\|_{\text{F}}^2 ~=~ (n^2-n)\cdot \frac{1}{(n-1)^2} ~=~ \frac{n}{n-1} ~=~ \max\limits_{\tr}(A).
\]

\item As all $3 \times 3$ doubly stochastic matrices saturating the inequality \eqref{E:EMR} are symmetric (up to permutations of the rows and columns), it is tempting to conjecture that this must be the case for any $n \times n$ doubly stochastic matrix saturating the inequality. However, this is not the case, as the following example shows. Let
\[
D := (1\oplus J_3) \times (J_2\oplus J_2) = \frac{1}{6}\begin{bmatrix} 3 & 3 & 0 & 0 \\ 1 & 1 & 2 & 2 \\ 1 & 1 & 2 & 2 \\ 1 & 1 & 2 & 2\end{bmatrix}
\]
One can readily see that no permutation of rows and/or columns will turn $D$ into a symmetric matrix. But
\[
\|D\|_{\text{F}}^2 ~=~ \frac{4}{3} ~=~ \max\limits_{\tr}(D).
\]

\item We saw that both the above matrix $D$ and the matrix
\[
\begin{bmatrix}
		0 & \frac{1}{2} & \frac{1}{2} \\[3pt]
		\frac{1}{2} & \frac{1}{4} & \frac{1}{4} \\[3pt]
		\frac{1}{2} & \frac{1}{4} & \frac{1}{4}
	\end{bmatrix}
\]
saturates the inequality \eqref{E:EMR}. Moreover, according to identity \eqref{E:matrix-prod}, both of them can be represented as the product of two matrices of the form
\begin{equation}\label{E:format-AB}
P(J_{n_1} \oplus \cdots \oplus J_{n_r})Q
\end{equation}
In fact, if $A$ and $B$ are two $n\times n$ doubly stochastic matrices of the form $P(J_{n_1} \oplus \cdots \oplus J_{n_r})Q$, where $P$ and $Q$ are permutation matrices and $n_1+\cdots+n_r=n$, then there exists a permutation matrix $R$ such that
\[
\|AB\|_{\operatorname{F}}^{2} = \tr(ABR).
\]
To verify this, let us write
\[
A=P_1 ( J_{n_1} \oplus \cdots \oplus J_{n_r})Q_1 =: P_1 A' Q_1,
\]
and
\[
B=P_2 ( J_{m_1} \oplus \cdots \oplus J_{m_r})Q_2 =: P_2 B' Q_2.
\]
Then $AB \,=\, P_1 A' Q_1 P_2 B' Q_2$. Hence, by the cyclic invariance of the trace, by the orthogonality of the permutation matrices, and by the fact that $A'$ and $B'$ are idempotent and symmetric, we see that
\begin{align*}
	\|AB\|_{\operatorname{F}}^2 ~&=~ \tr\left(AB(AB)^*\right) \\
	&=~ \tr\left( P_1 A' Q_1 P_2 B' Q_2 Q_2^\intercal  B'^\intercal  P_2^\intercal  Q_1^\intercal  A'^\intercal  P_1^\intercal  \right) \\
	&=~ \tr\left( A'^\intercal P_1^\intercal P_1 A' Q_1 P_2 B' B'^\intercal  P_2^\intercal  Q_1^\intercal  \right) \\
	&=~ \tr\left( A'^\intercal A' Q_1 P_2 B' P_2^\intercal  Q_1^\intercal  \right) \\
	&=~ \tr\left( A' Q_1 P_2 B' P_2^\intercal  Q_1^\intercal  \right) \\
	&=~ \tr\left( P_1^\intercal  P_1 A' Q_1 P_2 B' Q_2 Q_2^\intercal  P_2^\intercal  Q_1^\intercal  \right) \\
	&=~ \tr\left( AB Q_2^\intercal  P_2^\intercal  Q_1^\intercal  P_1^\intercal  \right) \\
	&=~ \tr\left( AB P \right),
\end{align*}
where $P:= \left(P_1 Q_1 P_2 Q_2\right)^\intercal $. However, these products do not always give matrices which saturate the inequality \eqref{E:EMR}. For example, consider the matrices $A:=J_3\oplus J_3 \oplus J_3$ and $B:= J_2\oplus J_3 \oplus J_4$. Then
\[
AB ~=~ \begin{bmatrix}
	1/3 & 1/3 & 1/9 & 1/9 & 1/9 & 0 & 0 & 0 & 0 \\
	1/3 & 1/3 & 1/9 & 1/9 & 1/9 & 0 & 0 & 0 & 0 \\
	1/3 & 1/3 & 1/9 & 1/9 & 1/9 & 0 & 0 & 0 & 0 \\
	0 & 0 & 2/9 & 2/9 & 2/9 & 1/12 & 1/12 & 1/12 & 1/12 \\
	0 & 0 & 2/9 & 2/9 & 2/9 & 1/12 & 1/12 & 1/12 & 1/12 \\
	0 & 0 & 2/9 & 2/9 & 2/9 & 1/12 & 1/12 & 1/12 & 1/12 \\
	0 & 0 & 0 & 0 & 0 & 1/4 & 1/4 & 1/4 & 1/4 \\
	0 & 0 & 0 & 0 & 0 & 1/4 & 1/4 & 1/4 & 1/4 \\
	0 & 0 & 0 & 0 & 0 & 1/4 & 1/4 & 1/4 & 1/4
\end{bmatrix}.
\]
It is easily verified that $\|AB\|_{\operatorname{F}}^2=\frac{37}{18}$ while $\max\limits_{\tr}(AB)=\frac{25}{12}$.%, which is obtained by permuting the third and sixth rows of $AB$.

Moreover, notice that not all doubly stochastic matrices which saturates the inequality \eqref{E:EMR} can be represented as the product of two matrices of the form $P(J_{n_1} \oplus \cdots \oplus J_{n_r})Q$. For example, it is easy to show that the matrix
\[
\begin{bmatrix}
	\tfrac{3}{5}&0&\tfrac{2}{5}\\[3pt]
	0&\tfrac{3}{5}&\tfrac{2}{5}\\[3pt]
	\tfrac{2}{5}&\tfrac{2}{5}&\tfrac{1}{5}
\end{bmatrix}
\]
is not of this form. In the light of these observations, we ask the following question:
\begin{question}
If $A$ and $B$ are of the form \eqref{E:format-AB}, when does $AB$ saturate the inequality \eqref{E:EMR}?
\end{question}

Finally, let us highlight that every known doubly stochastic solution to the equation $\|A\|_{\text{F}}^2 = \max\limits_{\tr}(A)$ has only rational coefficients. Hence, we propose the following question:
\begin{question}
Are the coefficients of every doubly stochastic solution of the equation $\|A\|_{\operatorname{F}}^2 = \smash{\max\limits_{\tr}}(A)$ rational?
\end{question}

\end{enumerate}

%%\begin{thebibliography}
%\bibliographystyle{plain}
%\bibliography{OnAQuestionOfErdos}

\begin{thebibliography}{10}
	
	\bibitem{MR540958}
	K.~Balasubramanian.
	\newblock Maximal diagonal sums.
	\newblock {\em Linear and Multilinear Algebra}, 7(3):249--251, 1979.
	
	\bibitem{Birkhoff1946}
	Garrett Birkhoff.
	\newblock Three observations on linear algebra.
	\newblock {\em Univ. Nac. Tucum\'{a}n. Revista A.}, 5:147--151, 1946.
	
	\bibitem{MR602332}
	G.~P. Egory\v{c}ev.
	\newblock {\em Reshenie problemy van-der-Vardena dlya permanentov}.
	\newblock Akad. Nauk SSSR Sibirsk. Otdel., Inst. Fiz., Krasnoyarsk, 1980.
	\newblock Preprint IFSO-13 M.
	
	\bibitem{MR642395}
	G.~P. Egory\v{c}ev.
	\newblock The solution of van der {W}aerden's problem for permanents.
	\newblock {\em Adv. in Math.}, 42(3):299--305, 1981.
	
	\bibitem{MR625097}
	D.~I. Falikman.
	\newblock Proof of the van der {W}aerden conjecture on the permanent of a
	doubly stochastic matrix.
	\newblock {\em Mat. Zametki}, 29(6):931--938, 957, 1981.
	
	\bibitem{MR117243}
	M.~Marcus and R.~Ree.
	\newblock Diagonals of doubly stochastic matrices.
	\newblock {\em Quart. J. Math. Oxford Ser. (2)}, 10:296--302, 1959.
	
	\bibitem{MR104679}
	Marvin Marcus and Morris Newman.
	\newblock On the minimum of the permanent of a doubly stochastic matrix.
	\newblock {\em Duke Math. J.}, 26:61--72, 1959.
	
	\bibitem{van1926aufgabe}
	Bartel~Leendert van~der Waerden.
	\newblock Aufgabe 45.
	\newblock {\em Jber. Deutsch. Math. Verein}, 35(117):23, 1926.
	
	\bibitem{MR672919}
	J.~H. van Lint.
	\newblock The van der {W}aerden conjecture: two proofs in one year.
	\newblock {\em Math. Intelligencer}, 4(2):72--77, 1982.
	
	\bibitem{Wang1974}
	Edward Tzu-Hsia Wang.
	\newblock Maximum and minimum diagonal sums of doubly stochastic matrices.
	\newblock {\em Linear Algebra Appl.}, 8(6):483--505, 1974.
	
\end{thebibliography}
%%\nocite{*}
%%\end{thebibliography}

\end{document}